\documentclass[11pt]{article}
\usepackage{mathrsfs}
\usepackage{amssymb}
\usepackage{amsmath}
\usepackage[all]{xy}
\def\Im{\mathop{\rm Im}\nolimits}
\def\Ker{\mathop{\rm Ker}\nolimits}
\def\Coker{\mathop{\rm Coker}\nolimits}
\def\Tr{\mathop{\rm Tr}\nolimits}
\def\mod{\mathop{\rm mod}\nolimits}
\def\grade{\mathop{\rm grade}\nolimits}
\def\rgrade{\mathop{\rm r.grade}\nolimits}
\def\fd{\mathop{\rm fd}\nolimits}
\def\id{\mathop{\rm id}\nolimits}
\def\pd{\mathop{\rm pd}\nolimits}
\def\gldim{\mathop{\rm gl.dim}\nolimits}
\def\add{\mathop{\rm add}\nolimits}
\def\Proj{\mathop{\rm Proj}\nolimits}
\def\cx{\mathop{\rm cx}\nolimits}
\def\inf{\mathop{\rm inf}\nolimits}
\def\sup{\mathop{\rm sup}\nolimits}
\def\dim{\mathop{\rm dim}\nolimits}

\title{\Large \bf The Existence of Maximal $n$-Orthogonal
Subcategories
}
\author{Zhaoyong Huang\thanks{{\it E-mail address}:
huangzy@nju.edu.cn}, Xiaojin Zhang\thanks{{\it E-mail address}:
xiaojinzhang@sohu.com}\\
{\footnotesize \it Department of Mathematics, Nanjing University,
Nanjing 210093, Jiangsu Province, P. R. China}}
\date{ }
\begin{document}
\baselineskip=18pt \maketitle
\begin{abstract}
For an $(n-1)$-Auslander algebra $\Lambda$ with global dimension
$n$, we give some necessary conditions for $\Lambda$ admitting a
maximal $(n-1)$-orthogonal subcategory in terms of the properties of
simple $\Lambda$-modules with projective dimension $n-1$ or $n$. For
an almost hereditary algebra $\Lambda$ with global dimension 2, we
prove that $\Lambda$ admits a maximal 1-orthogonal subcategory if
and only if for any non-projective indecomposable $\Lambda$-module
$M$, $M$ is injective is equivalent to that the reduced grade of $M$
is equal to 2. We give a connection between the Gorenstein Symmetric
Conjecture and the existence of maximal $n$-orthogonal subcategories
of $^{\bot}T$ for a cotilting module $T$. For a Gorenstein algebra,
we prove that all non-projective direct summands of a maximal
$n$-orthogonal module are $\Omega ^n\tau$-periodic. In addition, we
study the relation between the complexity of modules and the
existence of maximal $n$-orthogonal subcategories for the tensor
product of two finite-dimensional algebras.
\end{abstract}

\vspace{0.5cm}

\centerline{\large \bf 1. Introduction}

In [Iy2] and [Iy3], Iyama developed the classical 2-dimensional
Auslander-Reiten theory to higher-dimensional Auslander-Reiten
theory. For example, in [Iy2], the notion of almost split sequences
was generalized to that of $n$-almost split sequences; and in [Iy3],
the famous 2-dimensional Auslander correspondence was generalized to
higher-dimensional Auslander correspondence. In particular, in
[Iy2], Iyama introduced the notion of maximal $n$-orthogonal
subcategories, which played a crucial role in these two papers
mentioned above. In fact, Iyama's higher-dimensional
Auslander-Reiten theory depends on the existence of maximal
$n$-orthogonal subcategories. So, a natural question is: When do
maximal $n$-orthogonal subcategories exist? Geiss, Leclerc and
Schr\"oer proved in [GLS] that maximal $1$-orthogonal subcategories
exist for preprojective algebras as well as for certain algebras of
finite representation type. Erdmann and Holm gave in [EH] a
necessary condition that a selfinjective algebra admits a maximal
$n$-orthogonal subcategory. They proved that for a selfinjective
finite-dimensional $K$-algebra $\Lambda$, the complexity of any
$\Lambda$-module $M$ is at most 1 if $\Lambda$ admits a maximal
$n$-orthogonal module for some $n\geq 1$, where the {\it complexity}
of $M$, denoted by ${\rm cx}(M)$, is defined as ${\rm inf}\{ b \in
\mathbb{N}_0\ |$ there exists a $c>0$ such that ${\rm dim}_K
P_{n}\leq cn^{b-1}$ for all $n\}$ with $P_{n}$ the $(n+1)$st term in
a minimal projective resolution of $M$. This result means that for a
selfinjective algebra, maximal $n$-orthogonal modules rarely exist.
Recently, Iyama proved in [Iy4] that if $\Lambda$ is a
finite-dimensional algebra of finite representation type with
Auslander algebra $\Gamma$ and $\mod \Gamma$ contains a maximal
1-orthogonal object, then $\Lambda$ is hereditary and the dominant
dimension of $\Lambda$ is at least 1; furthermore, if the base field
is algebraically closed, then such $\Lambda$ is an upper triangular
matrix ring. The aim of the paper is to study algebras of global
dimension $n$ which admit maximal $(n-1)$-orthogonal subcategories
or modules, and that they are natural generalizations of Gabriel's
classification theorem on path algebras of Dynkin quivers (for
$n=1$). This paper is organized as follows.

In Section 2, we give some notions and notations and collect some
preliminary results.

In Section 3, we study the existence of maximal $n$-orthogonal
subcategories for algebras with finite global dimension. We observe
that an Artinian algebra $\Lambda$ with global dimension $n(\geq 1)$
admits no maximal $n$-orthogonal subcategories of $\mod \Lambda$.
Let $\Lambda$ be an $(n-1)$-Auslander algebra with global dimension
$n$. We prove that if $\Lambda$ admits a maximal $(n-1)$-orthogonal
subcategory of $\mod \Lambda$, then we can classify the simple
modules with projective dimension $n-1$; furthermore, in terms of
the properties of simple modules in $\mod \Lambda$ with projective
dimension $n$, we give some necessary and sufficient conditions for
$\Lambda$ admitting a trivial maximal $(n-1)$-orthogonal subcategory
of $\mod \Lambda$, and also give a necessary condition for $\Lambda$
admitting a non-trivial maximal $(n-1)$-orthogonal subcategory of
$\mod \Lambda$.

We also study the existence of maximal $n$-orthogonal subcategories
for algebras with global dimension 2. For example, if putting $n=1$,
then the above results holds true for classical Auslander algebras
with global dimension 2. In addition, for an almost hereditary
algebra $\Lambda$ with global dimension 2, we prove that $\Lambda$
admits a maximal 1-orthogonal subcategory of $\mod \Lambda$ if and
only if for any non-projective indecomposable module $M \in \mod
\Lambda$, $M$ is injective is equivalent to that the reduced grade
of $M$ is equal to 2. In particular, if $\Lambda$ admits a maximal
1-orthogonal subcategory $\mathscr{C}$ of $\mod \Lambda$, then
$\mathscr{C}=\add _{\Lambda}(\Lambda\bigoplus \mathbb{D}\Lambda
^{op})$.

In Section 4, we study the existence of maximal $n$-orthogonal
subcategories of $^{\bot}T$ for a cotilting module $T$. We give a
connection between the Gorenstein Symmetric Conjecture with the
existence of maximal $n$-orthogonal subcategories of $^{\bot}T$. As
a generalization of a result of Erdmann and Holm in [EH], we prove
that for a Gorenstein algebra, all non-projective direct summands of
a maximal $n$-orthogonal module are $\Omega ^n\tau$-periodic. It
should be pointed out that this result can be induced by [Iy3,
Theorem 2.5.1(1)].

In Section 5, both $\Lambda$ and $\Gamma$ are finite-dimensional
$K$-algebras over a field $K$. We prove that if the selfinjective
dimension of $\Lambda$ is equal to $n (\geq 1)$ and ${\rm
Hom}_{\Gamma}(\mathbb{D}\Gamma ^{op},\Gamma)\neq 0$, then ${\Lambda
\otimes_{K}\Gamma}$ admits no maximal $j$-orthogonal subcategories
of $\mod {\Lambda \otimes_{K}\Gamma}$ for any $j \geq n$. By this
result, we can construct algebras with infinite global dimension
admitting no maximal $n$-orthogonal subcategories for any $n \geq
1$. In addition, we prove that $\max\{\cx(M),\cx(N)\}\leq
\cx(M\otimes_{K}N)\leq \cx(M)+\cx(N)$ for any $M \in \mod \Lambda$
and $N \in \mod \Gamma$. As an application of this result, we can
construct a class of algebras $\Lambda$ with selfinjective dimension
$n(\geq 1)$, such that not all modules in $\mod \Lambda$ are of
complexity at most 1, but $\Lambda$ admits no maximal $j$-orthogonal
subcategories of $\mod \Lambda$ for any $j\geq n$.

\vspace{0.5cm}

\centerline{\large \bf 2. Preliminaries}

\vspace{0.2cm}

In this section, we give some notions and notations in our
terminology and collect some preliminary results for later use.

For a ring $\Lambda$, we use $\mod \Lambda$, $\gldim \Lambda$ and
$J(\Lambda)$ to denote the category of finitely generated left
$\Lambda$-modules, the global dimension and the Jacobson radical of
$\Lambda$, respectively. We use $\Tr$ and $(-)^*$ to denote the
Auslander transpose and the functor ${\rm Hom}_{\Lambda}(-,
\Lambda)$, respectively.

Let $M$ be a $\Lambda$-module. The {\it grade} of $M$, denoted by
$\grade M$, is defined as ${\rm inf}\{n\geq 0\ |\ {\rm
Ext}_{\Lambda}^{n}(M,\Lambda)\neq0\}$; and the {\it reduced grade}
of $M$, denoted by $\rgrade M$, is defined as ${\rm inf}\{n\geq 1\
|\ {\rm Ext}_{\Lambda}^{n}(M,\Lambda)\neq 0\}$ (see [Ho]). We use
$\pd_{\Lambda}M$, $\fd_{\Lambda}M$ and $\id _{\Lambda}M$ to denote
the projective, flat and injective dimensions of $M$, respectively.
Let $\Lambda$ be a left and right Noetherian ring and $M\in \mod
\Lambda$. We use
$$0\rightarrow M\rightarrow I^{0}(M)\rightarrow
I^{1}(M)\rightarrow \cdots \rightarrow I^{i}(M) \rightarrow \cdots$$
to denote a minimal injective resolution of $_{\Lambda}M$. For
positive integers $m$ and $n$, recall from [Iy1] that $\Lambda$ is
said to satisfy the {\it $(m,n)$-condition} (resp. the {\it
$(m,n)^{op}$-condition}) if $\fd _{\Lambda}I^i(\Lambda)$ (resp. $\fd
_{\Lambda^{op}}I^i(\Lambda ^{op})) \leq m-1$ for any $0 \leq i \leq
n-1$.

\vspace{0.2cm}

{\bf Lemma 2.1} ([Iy1, Proposition 2.4]) {\it Let $\Lambda$ be a
left and right Noetherian ring satisfying the $(n,n)$-condition and
$(n,n)^{op}$-condition. Then the subcategory $\{X \ |\ \grade X \geq
n\}$ {\it of} $\mod \Lambda$ is closed under submodules and factor
modules.}

\vspace{0.2cm}

Recall from [FGR] that a left and right Noetherian ring $\Lambda$ is
called {\it $n$-Gorenstein} if $\fd _{\Lambda}I^{i}(\Lambda)\leq i$
for any $0\leq i\leq n-1$. By [FGR, Theorem 3.7], the notion of
$n$-Gorenstein rings is left and right symmetric. It is clear that
$\Lambda$ is $n$-Gorenstein if and only if $\Lambda$ satisfies the
$(i, i)$ (or $(i, i)^{op}$)-condition for any $1 \leq i \leq n$.
Recall from [Bj] that $\Lambda$ is called {\it Auslander-Gorenstein}
if $\Lambda$ is $n$-Gorenstein for all $n$ and both $\id
_{\Lambda}\Lambda$ and $\id _{\Lambda ^{op}}\Lambda$ are finite; and
$\Lambda$ is called {\it Auslander-regular} if $\Lambda$ is
$n$-Gorenstein for all $n$ and $\gldim \Lambda$ is finite.

\vspace{0.2cm}

{\bf Lemma 2.2} ([IS, Corollary 7]) {\it Let $\Lambda$ be an
Auslander-Gorenstein ring with $\id _{\Lambda}\Lambda=\id
_{\Lambda^{op}}\Lambda =n\geq 1$. Then $\bigoplus _{i=0}^{n-1}
I^{i}(\Lambda)$ and $I^{n}(\Lambda)$ have no isomorphic direct
summands in common.}

\vspace{0.2cm}

The following easy observation is well known.

\vspace{0.2cm}

{\bf Lemma 2.3} {\it Let $\Lambda$ be a left Noetherian ring and
$M\in \mod \Lambda$ with $\pd _{\Lambda}M =n \ (<\infty)$. Then
${\rm Ext}_{\Lambda}^{n}(M, \Lambda) \neq 0$.}

\vspace{0.2cm}

{\bf Lemma 2.4} {\it Let $\Lambda$ be an Auslander-Gorenstein ring
with $\id _{\Lambda}\Lambda=\id _{\Lambda^{op}}\Lambda =n$. If $S
\in \mod \Lambda$ is simple with $\pd _{\Lambda}S=n$, then
$S\subseteq I^{n}(\Lambda)$ and $S\not\subseteq
I^{0}(\Lambda)\bigoplus\cdots\bigoplus I^{n-1}(\Lambda)$.}

\vspace{0.2cm}

{\it Proof.} Note that ${\rm Ext}^{i}_{\Lambda}(S,\Lambda)\cong {\rm
Hom}_{\Lambda}(S,I^{i}(\Lambda))$ for any $i\geq 0$. Because $\pd
_{\Lambda}S=n$, ${\rm Ext}^{n}_{\Lambda}(S,\Lambda)\neq 0$ by Lemma
2.3. So ${\rm Hom}_{\Lambda}(S,I^{n}(\Lambda))\neq 0$ and
$S\subseteq I^{n}(\Lambda)$. Then by Lemma 2.2, $S\not\subseteq
I^{0}(\Lambda)\bigoplus\cdots\bigoplus I^{n-1}(\Lambda)$.
$\hfill{\square}$

\vspace{0.2cm}

As a generalization of the notion of classical Auslander algebras,
Iyama introduced in [Iy2] the notion of $n$-Auslander algebras as
follows.

\vspace{0.2cm}

{\bf Definition 2.5} ([Iy4]) For a positive integer $n$, an Artinian
algebra $\Lambda$ is called an {\it $n$-Auslander algebra} if
$\gldim \Lambda\leq n+1$ and $I^{0}(\Lambda),I^{1}(\Lambda), \cdots,
I^{n}(\Lambda)$ are projective.

\vspace{0.2cm}

The notion of $n$-Auslander algebras is left and right symmetric by
[Iy4, Theorem 1.2]. It is trivial that $n$-Auslander algebras with
global dimension at most $n$ are semisimple. Note that the notion of
1-Auslander algebras is just that of classical Auslander algebras.
In the following, we assume that $n\geq 2$ when an $(n-1)$-Auslander
algebra is concerned.

\vspace{0.2cm}

{\bf Proposition 2.6} {\it Let $\Lambda$ be an $(n-1)$-Auslander
algebra with $\gldim \Lambda=n$. Then we have the following

(1) There does not exist an injective simple module $S \in \mod
\Lambda$ with $1\leq\pd _{\Lambda}S\leq n-1$. Dually, there does not
exist a projective simple module $S \in \mod \Lambda$ with $1\leq\id
_{\Lambda}S\leq n-1$.

(2) For a projective simple module $S\in \mod \Lambda$, ${\rm
Hom}_{\Lambda}(I^0(\Lambda), S)\neq 0$ if and only if $S$ is
injective; and ${\rm Hom}_{\Lambda}(I^0(\Lambda), S)=0$ if and only
if $\id _{\Lambda}S=n$.}

\vspace{0.2cm}

{\it Proof.} (1) Let $S \in \mod \Lambda$ be simple with $\pd
_{\Lambda}S=i$ with $1 \leq i \leq n-1$. Then by Lemma 2.3, we have
that ${\rm Hom}_{\Lambda}(S, I^{i}(\Lambda)) \cong {\rm
Ext}^{i}_{\Lambda}(S, \Lambda)\neq 0$. If $S$ is injective, then $S$
is isomorphic to a direct summand of $I^{i}(\Lambda)$ and hence $S$
is projective, which is a contradiction. Dually, we get the second
assertion.

(2) Let $S\in \mod \Lambda$ be a projective simple module. If $S$ is
injective, then we have an epimorphism ${\rm
Hom}_{\Lambda}(I^0(\Lambda), S) \rightarrow {\rm
Hom}_{\Lambda}(\Lambda, S)(\cong S) \rightarrow 0$ and so ${\rm
Hom}_{\Lambda}(I^0(\Lambda), S) \neq 0$. Conversely, if ${\rm
Hom}_{\Lambda}(I^0(\Lambda), S) \neq 0$, then $S$ is isomorphic to a
direct summand of $I^0(\Lambda)$ and hence $S$ is injective. The
second assertion follows from the first one and (1).
$\hfill{\square}$

\vspace{0.2cm}

{\bf Definition 2.7} ([AR1]) Let $\Lambda$ be an Artinian algebra.
Assume that $\mathscr{C} \subseteq \mathscr{D}$ are full
subcategories of $\mod \Lambda$ and $D\in \mod \Lambda$,
$C\in\mathscr{C}$. The morphism $f: D\to C$ is said to be a {\it
left $\mathscr{C}$-approximation} of $D$ if ${\rm Hom}_{\Lambda}(C,
C')\to {\rm Hom}_{\Lambda}(D, C')\to 0$ is exact for any
$C'\in\mathscr{C}$. The morphism $f: D\to C$ is said to be {\it left
minimal} if an endomorphism $g: C\to C$ is an automorphism whenever
$f=gf$. The subcategory $\mathscr{C}$ is said to be {\it covariantly
finite} in $\mathscr{D}$ if every module in $\mathscr{D}$ has a left
$\mathscr{C}$-approximation. The notions of {\it (minimal) right
$\mathscr{C}$-approximations} and {\it contravariantly finite
subcategories} of $\mathscr{D}$ may be defined dually. The
subcategory $\mathscr{C}$ is said to be {\it functorially finite} in
$\mathscr{D}$ if it is both covariantly finite and contravariantly
finite in $\mathscr{D}$.

\vspace{0.2cm}

The following result is due to T. Wakamatsu.

\vspace{0.2cm}

{\bf Lemma 2.8} ([AR1, Lemma 1.3]) {\it Let $\Lambda$ be an Artinian
algebra, $\mathscr{C}$ a full subcategory of $\mod \Lambda$ which is
closed under extensions and $D\in \mod \Lambda$. If $D
\stackrel{f}{\rightarrow} C \rightarrow Z \rightarrow 0$ is exact
with $f$ a minimal left $\mathscr{C}$-approximation of $D$, then
${\rm Ext}_{\Lambda}^{1}(Z, \mathscr{C})=0$.}

\vspace{0.5cm}

\centerline{\large\bf 3. Maximal $n$-orthogonal subcategories of
$\mod \Lambda$}

\vspace{0.2cm}

From now on, all algebras are Artinian algebras. For an Artinian
$R$-algebra $\Lambda$, we denote by $\mathbb{D}$ the ordinary
duality, that is, $\mathbb{D}(-)={\rm Hom}_R(-, I(R/J(R)))$, where
$I(R/J(R))$ is the injective envelope of $R/J(R)$.

In this section, we will mainly study the existence of maximal
$(n-1)$-orthogonal subcategories for algebras with global dimension
$n$, especially for $(n-1)$-Auslander algebras and for almost
hereditary algebras. In particular, we will give some necessary
conditions for $(n-1)$-Auslander algebras with global dimension $n$
admitting maximal $(n-1)$-orthogonal subcategories, and give some
necessary and sufficient conditions for almost hereditary algebras
with global dimension 2 admitting maximal 1-orthogonal
subcategories.

Let $\mathscr{C}$ be a full subcategory of $\mod \Lambda$ and $n$ a
positive integer. We denote $^{\bot_n}\mathscr{C}= \{ X\in \mod
\Lambda \ |\ {\rm Ext}_{\Lambda}^{i}(X, C)=0$ for any $C \in
\mathscr{C}$ and $1 \leq i \leq n \}$, and $\mathscr{C}^{\bot_n}= \{
X\in \mod \Lambda \ |\ {\rm Ext}_{\Lambda}^{i}(C, X)=0$ for any $C
\in \mathscr{C}$ and $1 \leq i \leq n \}$.

\vspace{0.2cm}

{\bf Definition 3.1} ([Iy2]) Let $\mathscr{C}\subseteq \mathscr{D}$
be full subcategories of $\mod \Lambda$ and $\mathscr{C}$
functorially finite in $\mathscr{D}$. For a positive integer $n$,
$\mathscr{C}$ is called a {\it maximal $n$-orthogonal subcategory}
of $\mathscr{D}$ if $\mathscr{C}={^{\bot_n}\mathscr{C}}\bigcap
\mathscr{D}=\mathscr{C}^{\bot_n}\bigcap \mathscr{D}$.

\vspace{0.2cm}

From the definition above, we get easily that both $\Lambda$ and
$\mathbb{D}\Lambda ^{op}$ are in any maximal $n$-orthogonal
subcategory of $\mod \Lambda$. For a module $M \in \mod \Lambda$, we
use $\add _{\Lambda}M$ to denote the subcategory of $\mod \Lambda$
consisting of all modules isomorphic to direct summands of finite
direct sums of copies of $_{\Lambda}M$.

\vspace{0.2cm}

{\bf Proposition 3.2} {\it Let $n$ be a positive integer. If $\id
_{\Lambda}\Lambda=n$ (in particular, if $\gldim \Lambda =n$), then
$\Lambda$ admits no maximal $n$-orthogonal subcategories of $\mod
\Lambda$.}

\vspace{0.2cm}

{\it Proof.} By the dual version of Lemma 2.3, we have that ${\rm
Ext}_{\Lambda}^{n}(\mathbb{D}\Lambda ^{op}, \Lambda) \neq 0$ and so
the assertion follows. $\hfill{\square}$

\vspace{0.2cm}

{\bf Proposition 3.3} {\it Let $\Lambda$ be an algebra with $\gldim
\Lambda =n\geq 2$. If $\Lambda$ admits a maximal $(n-1)$-orthogonal
subcategory of $\mod \Lambda$, then the following statements are
equivalent.

(1) $\Lambda$ is an $(n-1)$-Auslander algebra.

(2) $\pd _{\Lambda}\bigoplus_{i=1}^{n-1}I^i(\Lambda) \leq n-1$.

(2)$^{op}$ $\pd _{\Lambda ^{op}}\bigoplus_{i=1}^{n-1}I^i(\Lambda
^{op}) \leq n-1$.}

\vspace{0.2cm}

{\it Proof.} $(1) \Rightarrow (2)$ and $(1) \Rightarrow (2) ^{op}$
are trivial.

$(2) \Rightarrow (1)$ Because $\Lambda$ admits a maximal
$(n-1)$-orthogonal subcategory of $\mod \Lambda$, ${\rm
Ext}_{\Lambda}^{j}(\mathbb{D}\Lambda ^{op},\Lambda)=0$ for any
$1\leq j\leq n-1$. Then $\pd
_{\Lambda}\bigoplus_{i=1}^{n-1}I^i(\Lambda) \leq n-1$ implies
$\bigoplus_{i=1}^{n-1}I^i(\Lambda)$ is projective. Put $K=\Im
(I^0(\Lambda) \to I^1(\Lambda))$. Since $\gldim \Lambda =n$, $\pd
_{\Lambda}K \leq 1$. So $\pd _{\Lambda}I^0(\Lambda) \leq 1$ and
hence $I^0(\Lambda)$ is projective.

$(2)^{op}\Rightarrow (1)$  Note that $\mathscr{C}$ is a maximal
$(n-1)$-orthogonal subcategory of $\mod \Lambda$ if and only if
$\mathbb{D}\mathscr{C}$ is a maximal $(n-1)$-orthogonal subcategory
of $\mod \Lambda^{op}$. Then we get the assertion by using an
argument similar to that in the proof of $(2)\Rightarrow (1)$.
$\hfill{\square}$

\vspace{0.2cm}

Let $\Lambda$ be an algebra with $\gldim \Lambda =n\geq2$ admitting
a maximal $(n-1)$-orthogonal subcategory of $\mod \Lambda$. Then by
Proposition 3.3, we have that $\Lambda$ is an $(n-1)$-Auslander
algebra if and only if $\Lambda$ is Auslander-regular, if and only
if $\Lambda$ satisfies the $(n, n)$-condition, if and only if
$\Lambda$ satisfies the $(n, n)^{op}$-condition.

\vspace{0.2cm}

{\bf Proposition 3.4} {\it Let $\Lambda$ be an algebra with $\gldim
\Lambda=n\geq2$ and $\mathscr{C}$ a subcategory of $\mod \Lambda$
such that $\Lambda \in\mathscr{C}$ and ${\rm
Ext}_{\Lambda}^{j}(\mathscr{C},\mathscr{C})=0$ for any $1\leq j\leq
n-1$. Then $\grade M=n$ for any $M\in \mathscr{C}$ without
projective direct summands.}

\vspace{0.2cm}

{\it Proof.} Since $\Lambda \in\mathscr{C}$ and ${\rm
Ext}_{\Lambda}^{j}(\mathscr{C},\mathscr{C})=0$ for any $1 \leq j
\leq n-1$, $\mathscr{C} \subseteq {^{\bot _{n-1}}\Lambda}$. Notice
that $\gldim \Lambda=n$, so $M\in \mathscr{C}$ without projective
direct summands implies that $\pd _{\Lambda}M=n$. Let $0\rightarrow
P_{n}\rightarrow\cdots\rightarrow P_{1}\rightarrow P_{0}\rightarrow
M\rightarrow0$ be a minimal projective resolution of $M$. The
induced exact sequence $0\rightarrow M^{*}\rightarrow
P_{0}^{*}\rightarrow P_{1}^{*}\rightarrow\cdots\rightarrow
P_{n}^{*}\rightarrow \Tr \Omega^{n-1} M\rightarrow0$ with
$P_{0}^{*}\rightarrow P_{1}^{*}\rightarrow\cdots\rightarrow
P_{n}^{*}\rightarrow \Tr \Omega^{n-1} M\rightarrow0$ a minimal
projective resolution of $\Tr \Omega^{n-1} M$ by [M, Proposition
4.2], where $\Omega ^{n-1}M$ is the $(n-1)$st syzygy of $M$. Thus
$M^{*}=0$ by $\gldim \Lambda=n$ and therefore $\grade M=n$.
$\hfill{\square}$

\vspace{0.2cm}

The following is an immediate consequence of Proposition 3.4.

\vspace{0.2cm}

{\bf Corollary 3.5} {\it Let $\Lambda$ be an algebra with $\gldim
\Lambda =n\geq2$ admitting a maximal $(n-1)$-orthogonal subcategory
$\mathscr{C}$. Then $\grade M=n$ for any $M\in \mathscr{C}$ without
projective direct summands.}

\vspace{0.2cm}

{\bf Proposition 3.6} {\it Let $\Lambda$ be an Auslander-regular
algebra with $\gldim \Lambda=n$. Then there exist simple modules
$S_1, S_2 \in \mod \Lambda$ with $\pd _{\Lambda}S_1=n-1$ and $\pd
_{\Lambda}S_2=n$.}

\vspace{0.2cm}

{\it Proof.} If there do not exist simple modules with projective
dimension $n-1$, then we claim that ${\rm
Ext}_{\Lambda}^{n-1}(S,\Lambda)= 0$ for any simple module $S\in \mod
\Lambda$. If $\pd_{\Lambda}S\leq n-2$, then ${\rm
Ext}_{\Lambda}^{n-1}(S,\Lambda)= 0$. If $\pd _{\Lambda}S=n$, then
${\rm Hom}_{\Lambda}(S, I^{n-1}(\Lambda))=0$ by Lemma 2.4. So ${\rm
Ext}_{\Lambda}^{n-1}(S, \Lambda)=0$ and the claim is proved. It
follows that $\id _{\Lambda}\Lambda \leq n-2$, which is a
contradiction because $\gldim \Lambda =n$. The other assertion is
trivial. $\hfill{\square}$

\vspace{0.2cm}

Recall from [AR1] that a homomorphism $f: X \rightarrow Y$ in $\mod
\Lambda$ is called {\it left minimal} if an endomorphism $g: Y
\rightarrow Y$ is an automorphism whenever $f=gf$. Let $\mathscr{C}$
be a subcategory of $\mod \Lambda$. Recall from [Iy2] that a
complex:
$$M\stackrel{f_{0}}{\longrightarrow} C_{0}\stackrel{f_{1}}{\longrightarrow}
C_{1}\stackrel{f_{2}}{\longrightarrow}\cdots  \eqno{(1)}$$ with
$C_{i}\in\mathscr{C}$ for any $i\geq0$ is called {\it minimal} if
$f_i$ is left minimal for any $i \geq 0$. If the following sequence
of functors:
$$\cdots\stackrel{{\rm Hom}_{\Lambda}(f_2,\ )}{\longrightarrow}
{\rm Hom}_{\Lambda}(C_{1},\ )\stackrel{{\rm Hom}_{\Lambda}(f_1,\
)}{\longrightarrow} {\rm Hom}_{\Lambda}(C_{0},\ )\stackrel{{\rm
Hom}_{\Lambda}(f_0,\ )}{\longrightarrow} {\rm Hom}_{\Lambda}(M,\
)\rightarrow 0  $$ is exact on $\mathscr{C}$, then the complex (1)
is called a {\it left} $\mathscr{C}$-{\it resolution}. It is trivial
that if $\mathscr{C}$ is covariantly finite in $\mod \Lambda$, then
any $M\in \mod \Lambda$ has a minimal left $\mathscr{C}$-resolution.

 \vspace{0.2cm}

{\bf Lemma 3.7} {\it Let $\Lambda$ be an algebra with $\gldim
\Lambda=n$ admitting a maximal $(n-1)$-orthogonal subcategory
$\mathscr{C}$ of $\mod \Lambda$, and $X\in\mod \Lambda$. Then
$\id_{\Lambda} X\leq n-1$ if and only if the injective envelope
$I^0(X)$ of $X$ gives a minimal left $\mathscr{C}$-approximation of
$X$. In this case, the minimal injective resolution of $X$ gives a
minimal left $\mathscr{C}$-resolution of $X$.}

\vspace{0.2cm}

{\it Proof.} Let $X\in\mod \Lambda$ with $\id_{\Lambda} X=m\leq n-1$
and $C\in\mathscr{C}$. Since ${\rm
Ext}_{\Lambda}^{j}(\mathbb{D}\Lambda ^{op},C)=0$ for any $1\leq
j\leq n-1$, by applying the functor ${\rm Hom}_{\Lambda}(\ ,C)$ to a
minimal injective resolution of $X$, we get the following exact
sequence:
$$0\rightarrow {\rm Hom}_{\Lambda}(I^{m}(X), C)\rightarrow\cdots\rightarrow
{\rm Hom}_{\Lambda}(I^{1}(X), C)\rightarrow {\rm
Hom}_{\Lambda}(I^{0}(X), C) \rightarrow {\rm Hom}_{\Lambda}(X, C)
\rightarrow 0.$$ The necessity is proved.

Conversely, assume that the injective envelope $0\to X\to I^0(X)\to
Y\to0$ of $X$ gives a minimal left $\mathscr{C}$-approximation of
$X$. Since $\gldim \Lambda=n$, $\id_{\Lambda} Y\leq n-1$. By
assumption, the minimal injective resolution $0\to Y\to
I^0(Y)\to\cdots\to I^{n-1}(Y)\to 0$ of $Y$ gives a minimal left
$\mathscr{C}$-resolution of $Y$. Then we have a minimal left
$\mathscr{C}$-resolution $0\to X\to I^0(X)\to I^0(Y)\to\cdots\to
I^{n-1}(Y)\to 0$ of $X$. Since any module in $\mod \Lambda$ has a
minimal left $\mathscr{C}$-resolution of length at most $n$ (cf.
[Iy2, Theorem 2.2.3]), we have that $I^{n-1}(Y)=0$ and
$\id_{\Lambda} X\leq n-1$. $\hfill{\square}$

\vspace{0.2cm}

Let $\Lambda$ be an $(n-1)$-Auslander algebra with $\gldim
\Lambda=n$. Then by Proposition 3.6, there exist simple
$\Lambda$-modules with projective dimension $n-1$. On the other hand
hand, by Proposition 2.6, there do not exist injective simple
$\Lambda$-modules with projective dimension $n-1$. In terms of the
injective dimension of simple $\Lambda$-modules with projective
dimension $n-1$, we give some necessary conditions for $\Lambda$
admitting a maximal $(n-1)$-orthogonal subcategory of $\mod \Lambda$
as follows.

 \vspace{0.2cm}

{\bf Theorem 3.8} {\it Let $\Lambda$ be an $(n-1)$-Auslander algebra
with $\gldim \Lambda =n$ admitting a maximal $(n-1)$-orthogonal
subcategory $\mathscr{C}$ of $\mod \Lambda$ and $S \in \mod \Lambda$
simple with $\pd _{\Lambda}S=n-1$. Then we have

(1) $_{\Lambda}S$ is not injective.

(2) $1\leq\id _{\Lambda}S\leq n-1$ if and only if ${\rm
Hom}_{\Lambda}(S, P)=0$ for any indecomposable projective module in
$\mod \Lambda$ with $\id _{\Lambda}P=n$.

(3) $\id _{\Lambda}S=n$ if and only if ${\rm Hom}_{\Lambda}(S,P)\neq
0$ for some indecomposable projective module in $\mod \Lambda$ with
$\id _{\Lambda}P=n$.}

\vspace{0.2cm}

{\it Proof.} (1) Because $\pd _{\Lambda}S=n-1$, ${\rm
Ext}_{\Lambda}^{n-1}(S, \Lambda) \neq 0$ by Lemma 2.3. Notice that
$\Lambda$ admits a maximal $(n-1)$-orthogonal subcategory of $\mod
\Lambda$, so $S$ is not injective.

(2) We first prove the sufficiency. Let $S \in \mod \Lambda$ be a
simple module with $\pd _{\Lambda}S=n-1$ and $N$ a module in
$\mathscr{C}$ without projective direct summands. Then $\grade N=n$
by Corollary 3.5. We claim that ${\rm Hom}_{\Lambda}(S, N)=0$.
Otherwise, if ${\rm Hom}_{\Lambda}(S, N) \neq 0$, then $S$ is
isomorphic to a submodule of $N$ and so $\grade S=n$ by Lemma 2.1.
It follows that ${\rm Ext}_{\Lambda}^{n-1}(S, \Lambda)=0$, which is
a contradiction because $\pd _{\Lambda}S=n-1$. The claim is proved.
On the other hand, ${\rm Hom}_{\Lambda}(S, P)=0$ for any
indecomposable projective module in $\mod \Lambda$ with $\id
_{\Lambda}P=n$ by assumption.

Let $0 \to S \to C$ be a minimal left $\mathscr{C}$-approximation of
$S$ and $N$ a non-zero indecomposable direct summand of $C$. Then
${\rm Hom}_{\Lambda}(S, N)\neq 0$. By the above argument, $N$ is
projective and $\id _{\Lambda}N<n$. Because $\mathscr{C}$ is a
maximal $(n-1)$-orthogonal subcategory of $\mod \Lambda$ by
assumption, $N$ is injective. So $N\cong I^0(S)$ and hence $C\cong
I^0(S)$ by the minimality of the above left
$\mathscr{C}$-approximation of $S$. Then it follows from (1) and
Lemma 3.7 that $1\leq\id _{\Lambda}S\leq n-1$.

We next prove the necessity. Let $S \in \mod \Lambda$ be a simple
module with $\pd _{\Lambda}S=n-1$ and $P \in \mod \Lambda$ an
indecomposable projective module with $\id _{\Lambda}P=n$. If ${\rm
Hom}_{\Lambda}(S,P)\neq 0$, then there exists an epimorphism
$\mathbb{D}P \rightarrow \mathbb{D}S \rightarrow 0$. Since
$\mathscr{C}$ is a maximal $(n-1)$-orthogonal subcategory of $\mod
\Lambda$, $\mathbb{D}\mathscr{C}$ is a maximal $(n-1)$-orthogonal
subcategory of $\mod \Lambda ^{op}$. So $\grade \mathbb{D}P=n$ by
the opposite version of Proposition 3.4 and hence $\grade
\mathbb{D}S=n$ by the opposite version of Lemma 2.1. Since $\gldim
\Lambda =n$, $\pd _{\Lambda ^{op}}\mathbb{D}S=n$. So $\id
_{\Lambda}S=n$, which is a contradiction.

(3) Since $\Lambda$ is an $(n-1)$-Auslander algebra, the injective
dimension of simple modules in $\mod \Lambda$ with projective
dimension $n-1$ is situated between 1 and $n$ by Proposition 2.6.
Thus the assertion follows immediately from (2). $\hfill{\square}$

\vspace{0.2cm}

In the following, we will study the properties of simple modules in
$\mod \Lambda$ with $\pd _{\Lambda}S=n$ if an $(n-1)$-Auslander
algebra $\Lambda$ with $\gldim \Lambda=n$ admits a maximal
$(n-1)$-orthogonal subcategory of $\mod \Lambda$.

\vspace{0.2cm}

{\bf Proposition 3.9} {\it Let $\Lambda$ be an $(n-1)$-Auslander
algebra with $\gldim \Lambda =n$ admitting a maximal
$(n-1)$-orthogonal subcategory $\mathscr{C}$ of $\mod \Lambda$ and
$S$ a simple $\Lambda$-module with $\pd _{\Lambda}S=n$. Then we have

(1) $S$ is injective if and only if $S\in \mathscr{C}$ and ${\rm
Hom}_{\Lambda}(S,C)=0$ for any (non-projective) indecomposable
module $C \in \mathscr{C}$ with $C\not\cong I^0(S)$.

(2) $1\leq\id _{\Lambda}S\leq n-1$ if and only if $S\not\in
\mathscr{C}$ and ${\rm Hom}_{\Lambda}(S,C)=0$ for any
(non-projective) indecomposable module $C \in \mathscr{C}$ with
$C\not\cong I^0(S)$.

(3) $\id _{\Lambda}S=n$ if and only if ${\rm
Hom}_{\Lambda}(S,C)\neq0$ for some (non-projective) indecomposable
module $C \in \mathscr{C}$ with $C\not\cong I^0(S)$.}

\vspace{0.2cm}

{\it Proof.} (1) The necessity is easy. Conversely, assume that
$S\in \mathscr{C}$ and ${\rm Hom}_{\Lambda}(S,C)=0$ for any
non-projective indecomposable module $C \in \mathscr{C}$ with
$C\not\cong I^0(S)$. Because ${\rm Hom}_{\Lambda}(S,S)\neq 0$ and
$\pd _{\Lambda}S=n$ by assumption, $S\cong I^0(S)$ is injective.

(2) If $1\leq\id _{\Lambda}S\leq n-1$, then $S\not \in \mathscr{C}$.
Since $\mathscr{C}$ is maximal $(n-1)$-orthogonal, the minimal
injective resolution of $S$: $0\rightarrow
S\stackrel{f}{\rightarrow} I^0(S)\rightarrow I^0(S)/S\rightarrow 0$
is a minimal left $\mathscr{C}$-approximation of $S$ by Lemma 3.7.
So, for any indecomposable module $C\in \mathscr{C}$ with ${\rm
Hom}_{\Lambda}(S,C)\neq 0$ and $g\in {\rm Hom}_{\Lambda}(S, C)$,
there exists an $h\in {\rm Hom}_{\Lambda}(I^0(S),C)$ such that
$g=hf$. Since $f$ is essential and $g$ is a monomorphism, $h$ is a
splittable monomorphism and $C\cong I^0(S)$ is injective.

Conversely, since $\pd _{\Lambda}S=n$, $\grade S=n$ by Lemmas 2.3
and 2.4. So $S$ cannot be embedded into any projective module. On
the other hand, ${\rm Hom}_{\Lambda}(S,C)=0$ for any non-projective
indecomposable module $C \in \mathscr{C}$ with $C\not\cong I^0(S)$
by assumption. Since $S\not\in\mathscr{C}$, the injective envelope
$0\to S\rightarrow I^0(S)$ of $S$ is a minimal left
$\mathscr{C}$-approximation of $S$ by the argument of Theorem
3.8(2). Then by Lemma 3.7, the assertion follows.

(3) It follows from (1) and (2).  $\hfill{\square}$

\vspace{0.2cm}

For a positive integer $n$, we know that
$\add_{\Lambda}(\Lambda\bigoplus \mathbb{D}\Lambda ^{op})$ is
contained in any maximal $n$-orthogonal subcategory of $\mod
\Lambda$. On the other hand, it is easy to see that if
$\add_{\Lambda}(\Lambda\bigoplus \mathbb{D}\Lambda ^{op})$ is a
maximal $n$-orthogonal subcategory of $\mod \Lambda$, then
$\add_{\Lambda}(\Lambda\bigoplus \mathbb{D}\Lambda ^{op})$ is the
unique maximal $n$-orthogonal subcategory of $\mod \Lambda$. In this
case, we say that $\Lambda$ admits a {\it trivial maximal
$n$-orthogonal subcategory} of $\mod \Lambda$. As an application of
Proposition 3.9, we give some necessary and sufficient conditions
that an $(n-1)$-Auslander algebra $\Lambda$ with $\gldim \Lambda =n$
admits a trivial maximal $(n-1)$-orthogonal subcategory of $\mod
\Lambda$.

\vspace{0.2cm}

{\bf Corollary 3.10} {\it Let $\Lambda$ be an $(n-1)$-Auslander
algebra with $\gldim \Lambda=n$. Then the following statements are
equivalent.

(1) $\Lambda$ admits a trivial maximal $(n-1)$-orthogonal
subcategory $\add_{\Lambda}(\Lambda\bigoplus \mathbb{D}\Lambda
^{op})$ of $\mod \Lambda$.

(2) A simple module $S \in \mod \Lambda$ is injective if $\pd
_{\Lambda}S=n$.

(3) There do not exist simple modules in $\mod \Lambda$ with both
projective and injective dimensions $n$; and $1\leq\pd
_{\Lambda}S\leq n-1$ if and only if $1\leq\id _{\Lambda}S\leq n-1$
for a simple module $S \in \mod \Lambda$.}

\vspace{0.2cm}

{\it Proof.} $(1) \Rightarrow (2)$ Let $S \in \mod \Lambda$ be a
simple module with $\pd _{\Lambda}S=n$. Because
$\add_{\Lambda}(\Lambda\bigoplus \mathbb{D}\Lambda ^{op})$ is a
maximal $(n-1)$-orthogonal subcategory of $\mod \Lambda$, ${\rm
Hom}_{\Lambda}(S,C)=0$ for any non-projective indecomposable module
$C \in \mathscr{C}$ with $C\not\cong I^0(S)$. Then it follows from
Proposition 3.9 that $\id _{\Lambda}S \leq n-1$. On the other hand,
$\Lambda$ is an $(n-1)$-Auslander algebra, so ${\rm
Ext}_{\Lambda}^{j}(S,\Lambda)=0$ for any $1\leq j\leq n-1$ by Lemmas
2.3 and 2.4 and hence $S \in \mathscr{C}$. It follows that $S$ is
injective.

$(2) \Rightarrow (1)$ Let $E$ be an indecomposable direct summand of
$I^n(\Lambda)$. Then $E \cong I^0(S)$ for some simple module $S \in
\mod \Lambda$. So ${\rm Ext}_{\Lambda}^{n}(S, \Lambda)\cong {\rm
Hom}_{\Lambda}(S,I^n(\Lambda))\neq 0$ and $\pd _{\Lambda}S=n$. Thus
by (2), $S$ is injective and $E\cong S$. Then it follows easily from
Lemma 2.4 that $\grade E=n$, which implies that $\grade I^n(\Lambda)
=n$. On the other hand, because $\Lambda$ is an $(n-1)$-Auslander
algebra, $I^0(\Lambda), \cdots, I^{n-1}(\Lambda)$ are projective.
Notice that $\mathbb{D}\Lambda ^{op} \in \add _{\Lambda}\bigoplus
_{i=0}^{n}I^{i}(\Lambda)$, so ${\rm
Ext}_{\Lambda}^{j}(\mathbb{D}\Lambda ^{op},\Lambda)=0$ for any
$1\leq j\leq n-1$.

Now let $M\in {^{\bot_{n-1}}\Lambda}$ be indecomposable. If $\pd
_{\Lambda}M\leq n-1$, then $M$ is projective. If $\pd
_{\Lambda}M=n$, then $\grade M=n$. By Lemma 2.1, for a simple
submodule $S$ of $M$, $\grade S=n$. So $\pd _{\Lambda}S=n$ and hence
$S$ is injective by (2). It follows that $M (\cong S)$ is injective.
Thus we have that $\add_{\Lambda}(\Lambda\bigoplus \mathbb{D}\Lambda
^{op})={^{\bot_{n-1}}\Lambda}$ and therefore
$\add_{\Lambda}(\Lambda\bigoplus \mathbb{D}\Lambda ^{op})$ is a
maximal $(n-1)$-orthogonal subcategory of $\mod \Lambda$ by [Iy2,
Proposition 2.2.2].

$(3) \Rightarrow (2)$ It is easy.

$(1)+(2) \Rightarrow (3)$ It suffices to prove the latter assertion
by (2). If $S\in \mod \Lambda$ is simple with $1\leq\pd
_{\Lambda}S\leq n-1$, then we only have to show $\id _{\Lambda}S\neq
n$ by Proposition 2.6. Otherwise, if $\id _{\Lambda}S=n$, then $\pd
_{\Lambda ^{op}}\mathbb{D}S=n$. By (1), $\Lambda ^{op}$ admits a
trivial maximal $(n-1)$-orthogonal subcategory of $\mod \Lambda
^{op}$, so the opposite version of (2) holds true. Then it follows
that $\mathbb{D}S$ is injective and $S$ is projective, which is a
contradiction. The converse can be proved dually. $\hfill{\square}$

\vspace{0.2cm}

We give an example to illustrate Theorem 3.8 and Corollary 3.10.

\vspace{0.2cm}

{\bf Example 3.11 } Let $\Lambda$ be a finite-dimensional algebra
given by
 the quiver $Q$:  $$\xymatrix{1 & \ar[l]_{\beta _{1}} 2 & \ar[l]_{\beta _{2}} 3
& \ar[l]_{\beta _{3}} \cdots & \ar[l]_{\beta _{n}} n+1}
$$ modulo the ideal generated by $\{
\beta_{i}\beta_{i+1}| 1\leq i\leq n-1 \}$. Then we have

(1) $\Lambda$ is an $(n-1)$-Auslander algebra with $\gldim
\Lambda=n$ and admits a maximal $(n-1)$-orthogonal subcategory
$\mathscr{C}=\add_{\Lambda}\big(P(1)\bigoplus P(2)\bigoplus
P(3)\bigoplus\cdots\bigoplus P(n+1)\bigoplus S(n+1)\big)$.

(2) $\pd_{\Lambda} S(n)= n-1$, $\id_{\Lambda} P(1)= n$, ${\rm
Hom}_{\Lambda}(S(n),P(1))=0$ and $\id_{\Lambda} S(n)= 1$.

(3) $\{S\in \mod \Lambda |\ S$ is simple with $1\leq \pd
_{\Lambda}S\leq n-1 \}=\{S\in \mod \Lambda |\ S$  is simple with
$1\leq \id_{\Lambda}S\leq n-1\}=\{S(i)|2\leq i\leq n\}$.

(4) $\pd _{\Lambda}S(n+1)=n$ and $S(n+1)=I(n+1)$.

\vspace{0.2cm}

As another application of Proposition 3.9, we give a necessary
condition for an $(n-1)$-Auslander algebra $\Lambda$ with $\gldim
\Lambda =n$ admitting a non-trivial maximal $(n-1)$-orthogonal
subcategory of $\mod \Lambda$.

\vspace{0.2cm}

{\bf Corollary 3.12} {\it Let $\Lambda$ be an $(n-1)$-Auslander
algebra with $\gldim \Lambda=n$. If $\Lambda$ admits a non-trivial
maximal $(n-1)$-orthogonal subcategory $\mathscr{C}(\neq
\add_{\Lambda}(\Lambda\bigoplus \mathbb{D}\Lambda ^{op}))$ of $\mod
\Lambda$, then there exists a simple module $S\in \mod \Lambda$ such
that $\pd _{\Lambda}S=n$ and $\id _{\Lambda}S=n$.}

\vspace{0.2cm}

{\it Proof.} Let $\mathscr{C}\neq \add_{\Lambda}(\Lambda\bigoplus
\mathbb{D}\Lambda ^{op})$ be a maximal $(n-1)$-orthogonal
subcategory of $\mod \Lambda$. Then there exists an indecomposable
module $M\in\mathscr{C}$ such that $\pd _{\Lambda}M=n$ and $\id
_{\Lambda}M=n$. So $\grade M=n$. Then by Lemma 2.1, for any simple
submodule $S$ of $M$, $\grade S=n$ and $\pd _{\Lambda}S=n$. Notice
that $M\not \cong I^0(S)$ and $M\in\mathscr{C}$, so $\id
_{\Lambda}S=n$ by Proposition 3.9. $\hfill{\square}$

\vspace{0.2cm}

In the following, we will study the existence of maximal
1-orthogonal subcategories for some kinds of algebras (Auslander
algebras and almost hereditary algebras) with global dimension 2.
First, if putting $n=1$, we get the following corollary immediately
from Theorem 3.8. This result means that the existence of maximal
1-orthogonal subcategories of $\mod \Lambda$ for an Auslander
algebra $\Lambda$ with $\gldim \Lambda =2$ enables us to classify
the simple modules in $\mod \Lambda$ with projective dimension 1.

\vspace{0.2cm}

{\bf Corollary 3.13} {\it Let $\Lambda$ be an Auslander algebra with
$\gldim \Lambda =2$ admitting a maximal 1-orthogonal subcategory
$\mathscr{C}$ of $\mod \Lambda$ and $S \in \mod \Lambda$ simple with
$\pd _{\Lambda}S=1$. Then we have

(1) $\id _{\Lambda}S=1$ if and only if ${\rm Hom}_{\Lambda}(S, P)=0$
for any indecomposable projective module in $\mod \Lambda$ with $\id
_{\Lambda}P=2$.

(2) $\id _{\Lambda}S=2$ if and only if ${\rm Hom}_{\Lambda}(S,P)\neq
0$ for some indecomposable projective module in $\mod \Lambda$ with
$\id _{\Lambda}P=2$.}

\vspace{0.2cm}

{\bf Definition 3.14} ([HRS]) An algebra $\Lambda$ is called {\it
almost hereditary} if the following conditions are satisfied: (1)
$\gldim \Lambda\leq 2$; and (2) If $X \in \mod \Lambda$ is
indecomposable, then either $\id {_{\Lambda}X}\leq 1$ or $\pd
{_{\Lambda}X}\leq 1$.

\vspace{0.2cm}

By Proposition 3.2, if $\gldim \Lambda =1$, then $\Lambda$ admits no
maximal $n$-orthogonal subcategories of $\mod \Lambda$ for any
$n\geq 1$. For an almost hereditary algebra $\Lambda$ with $\gldim
\Lambda =2$, we give some equivalent characterizations of the
existence of maximal 1-orthogonal subcategories of $\mod \Lambda$ as
follows.

\vspace{0.2cm}

{\bf Theorem 3.15} {\it Let $\Lambda$ be an almost hereditary
algebra with $\gldim \Lambda =2$. Then the following statements are
equivalent.

(1) $\Lambda$ admits a maximal 1-orthogonal subcategory
$\mathscr{C}$ of $\mod \Lambda$.

(2) The following conditions are satisfied:

\ \ \ \ \ (i) $\rgrade \mathbb{D}\Lambda ^{op} =2$; and

\ \ \ \ \ (ii) Any non-projective indecomposable module $M\in \mod
\Lambda$ is injective if $\rgrade M$

\ \ \ \ \ \ \ \ \ \ $=2$.

(3) For any non-projective indecomposable module $M\in \mod
\Lambda$, $M$ is injective if and only if $\rgrade M=2$.

In particular, if $\Lambda$ admits a maximal 1-orthogonal
subcategory $\mathscr{C}$ of $\mod \Lambda$, then $\mathscr{C}=\add
_{\Lambda}(\Lambda\bigoplus \mathbb{D}\Lambda ^{op})$. That is, the
maximal 1-orthogonal subcategory of $\mod \Lambda$ is trivial if it
exists.}

\vspace{0.2cm}

{\it Proof.} $(1)\Rightarrow (2)$ Let $\mathscr{C}$ be a maximal
1-orthogonal subcategory of $\mod \Lambda$ and $N \in \mathscr{C}$
an indecomposable module. Because $\Lambda$ is almost hereditary,
$N$ is projective or injective. So $\mathscr{C}=\add
_{\Lambda}(\Lambda\bigoplus \mathbb{D}\Lambda
^{op})={^{\bot_{1}}\mathscr{C}}={^{\bot_{1}}\Lambda}$. Since $\pd
_{\Lambda} \mathbb{D}\Lambda ^{op} =\gldim \Lambda =2$, ${\rm
Ext}_{\Lambda}^2(\mathbb{D}\Lambda ^{op} , \Lambda)\neq 0$ by Lemma
2.3. So $\rgrade \mathbb{D}\Lambda ^{op} =2$. If $M\in \mod \Lambda$
is indecomposable with $\rgrade M=2$, then $M \in \mathscr{C}$; so,
if $M$ is further non-projective, then $M$ is injective by the above
argument.

$(2)\Rightarrow (3)$ It is easy.

$(3)\Rightarrow (1)$ We will show that $\mathscr{C}=\add
_{\Lambda}(\Lambda\bigoplus \mathbb{D}\Lambda ^{op})$ is a maximal
1-orthogonal of $\mod \Lambda$. By [Iy2, Proposition 2.2.2], we only
need to show ${^{\bot_{1}}\Lambda} \subseteq \add
_{\Lambda}(\Lambda\bigoplus \mathbb{D}\Lambda ^{op})$. Let $M\in
{^{\bot_{1}}\Lambda}$ be an indecomposable module. If $\pd
{_{\Lambda}M}=2$, then $\rgrade M=2$ and so $M$ is injective by
assumption. If $\pd {_{\Lambda}M}\leq 1$, then $M$ is projective.
Thus we conclude that ${^{\bot_{1}}\Lambda} \subseteq \add
_{\Lambda}(\Lambda\bigoplus \mathbb{D}\Lambda ^{op})$.

The last assertion follows from the above argument.
$\hfill{\square}$

\vspace{0.2cm}

At the end of this section, we give two examples related to
Corollary 3.13 and Theorem 3.15 as follows.

\vspace{0.2cm}

{\bf Example 3.16} Putting $n=2$ in Example 3.11, then we have

(1) $\Lambda$ is an Auslander algebra and an almost hereditary
algebra with $\gldim \Lambda =2$ admitting a maximal 1-orthogonal
subcategory $\mathscr{C}=\add_{\Lambda}\big(P(1) \bigoplus
P(2)\bigoplus P(3)\bigoplus S(3)\big)$.

(2) $\pd _{\Lambda}S(2)=1=\id _{\Lambda}S(2)$, $\id
_{\Lambda}P(1)=2$ and ${\rm Hom}_{\Lambda}(S(2),P(1))=0$.

(3) $S(3)=I(3)$ and $\pd _{\Lambda}S(3)=2$.

(4) $\rgrade S(3)=2$.

\vspace{0.2cm}

{\bf Example 3.17} Let $\Lambda$ be a finite-dimensional algebra
given by the quiver $Q:$
 $$\xymatrix{&6 \ar[r]^{\alpha}\ar[d]^{\gamma} & 4\ar[d]^{\beta}& \\
&5\ar[r]^{\delta}&3\ar[r]^{\lambda}\ar[d]^{\mu}&1\\
&\ \ &2&}$$ modulo the ideal generated by $\{
\beta\alpha-\delta\gamma, \mu\delta, \lambda\beta\}$. Then we have

(1) $\Lambda$ is an Auslander algebra and an almost hereditary
algebra with $\gldim \Lambda =2$.

(2) $\pd _{\Lambda}S(3)=1$, $\id _{\Lambda}S(3)=2$, $P(4)$, $P(5)$
and $P(6)$ are injective; $\id _{\Lambda}P(i)=2$ and ${\rm
Hom}_{\Lambda}(S(3),P(i))=0$ for $i=1,2,3$.

(3) ${\rm Ext}_{\Lambda}^{1}(\mathbb{D}\Lambda ^{op},\Lambda)=0$.
Both $S(4)$ and $S(5)$ are of $\rgrade$ 2, and neither of them are
injective.

(4) There does not exist a simple module $S$ such that $\pd
_{\Lambda}S=2$ and $\id _{\Lambda}S=2$.

So there exist no maximal 1-orthogonal subcategories of $\mod
\Lambda$ by Corollary 3.13 or Theorem 3.15.

\vspace{0.5cm}

\centerline{\large\bf 4. Maximal $n$-orthogonal subcategories of
$^{\bot}T$}

\vspace{0.2cm}

In this section, we will study the properties of an algebra
$\Lambda$ and $\Lambda$-modules if $\Lambda$ admits maximal
$n$-orthogonal subcategories of $^{\bot}T$ for a cotilting module
$T$.

Let $\mathscr{C}$ be a full subcategory of $\mod \Lambda$. We use
$\widehat{\mathscr{C}}$ to denote the subcategory of $\mod \Lambda$
consisting of the module $X$ for which there exists an exact
sequence $0\rightarrow C_{n}\rightarrow
C_{n-1}\rightarrow\cdots\rightarrow C_{0}\rightarrow X\rightarrow 0$
with each $C_{i}$ in $\mathscr{C}$. Denote
$^{\bot}\mathscr{C}=\bigcap _{n \geq 1} {^{\bot_n}\mathscr{C}}$, and
$\mathcal{I}^{\infty}(\Lambda)=\{ X\in \mod \Lambda \ |\ \id
_{\Lambda}X <\infty \}$.

\vspace{0.2cm}

{\bf Definition 4.1} ([AR1]) A module $T\in {\rm mod}\ \Lambda$ is
called a {\it cotilting} module if the following conditions are
satisfied: (1) $\id _{\Lambda}T=n<\infty$; (2) $T\in {^{\bot}T}$;
and (3) $\mathbb{D}\Lambda ^{op}\in \widehat{\add_{\Lambda}T}$. A
cotilting module $T$ is called {\it strong cotilting} if
$\mathcal{I}^{\infty}(\Lambda)=\widehat{\add _{\Lambda}T}$.

\vspace{0.2cm}

{\bf Lemma 4.2} {\it Let $\Lambda$ be an algebra and $T \in \mod
\Lambda$ a cotilting module. Then for any $M\in {^{\bot}T}$ with
$\id _{\Lambda}M=n$, there exists an exact sequence $0\rightarrow
M\rightarrow T_{0}\rightarrow T_{1}\rightarrow\cdots\rightarrow
T_{n}\rightarrow0$ with $T_i \in \add _{\Lambda}T$ for any $0 \leq i
\leq n$.}

\vspace{0.2cm}

{\it Proof.} For any $M\in {^{\bot}T}$, we have an exact sequence
$0\rightarrow M\stackrel{f_{0}}{\longrightarrow}
T_{0}\stackrel{f_{1}}{\longrightarrow}
T_{1}\stackrel{f_{2}}{\longrightarrow}\cdots\stackrel{f_{n}}{\longrightarrow}
T_{n}\stackrel{f_{n+1}}{\longrightarrow}\cdots$ with $\Coker
f_{i}\in {^{\bot}T}$ for any $i\geq 0$ by [AR1, Theorem 5.4]. Then
for any $N\in {^{\bot}T}$ and $i\geq1$, $0={\rm
Ext}_{\Lambda}^{n+i}(N,M)={\rm Ext}_{\Lambda}^{i}(N,\Coker
f_{n-1})$. So $\Coker f_{n-1}\in {(^{\bot}T)^{\bot}}$. Since $T$ is
cotilting and $\Coker f_{n-1}\in {^{\bot}T}$, $\Coker f_{n-1}\in
\add_{\Lambda}T$. $\hfill{\square}$

\vspace{0.2cm}

Auslander and Reiten in [AR1, p.150] posed a question: whether does
$\id _{\Lambda}\Lambda <\infty$ imply $\id _{\Lambda ^{op}}\Lambda
<\infty$? This question is now referred to the Gorenstein Symmetry
Conjecture, which still remains open (see [BR]). The following
result gives a connection between this conjecture with the existence
of maximal $n$-orthogonal subcategories. For a module $M\in \mod
\Lambda$, the {\it basic submodule} of $M$, denoted by $M_b$, is
defined as the direct sum of one copy of each non-isomorphic
indecomposable direct summand of $M$.

\vspace{0.2cm}

{\bf Theorem 4.3} {\it Let $\Lambda$ be an algebra with $\id
_{\Lambda}\Lambda=n\geq 1$ and $T \in \mod \Lambda$ a cotilting
module. If $T\in {^{\bot _n} \Lambda}$ (in particular, if $\Lambda$
admits a maximal $j$-orthogonal subcategory of ${^{\bot}T}$ for some
$j\geq n$), then $\id _{\Lambda ^{op}}\Lambda =n$ and $T$ is a
strong cotilting module with $T_{b}=\Lambda_{b}$.}

\vspace{0.2cm}

{\it Proof.} Let $T\in {^{\bot _n} \Lambda}$ be a cotilting module.
By Lemma 4.2, there exists an exact sequence $0\rightarrow
\Lambda\rightarrow T_{0}\rightarrow
T_{1}\rightarrow\cdots\rightarrow T_{n}\rightarrow 0$ with $T_i \in
\add _{\Lambda}T$ for any $0 \leq i \leq n$, which is splittable. So
we have $_{\Lambda}\Lambda \in \add _{\Lambda}T$. Note that all
cotilting modules in $\mod \Lambda$ have the same number of
non-isomorphic indecomposable direct summands which is equal to the
number of non-isomorphic indecomposable projective modules in $\mod
\Lambda$ (see [AR1] or [M]). It follows that $_{\Lambda}\Lambda _b
={_{\Lambda}T_b}$ and $_{\Lambda}\Lambda$ is a cotilting module.
Thus $\id _{\Lambda ^{op}}\Lambda =n$ by [AR2, Lemma 1.7].
$\hfill{\square}$

\vspace{0.2cm}

From now on, $\Lambda$ is a Gorenstein algebra, that is, $\id
_{\Lambda}\Lambda=\id _{\Lambda ^{op}}\Lambda <\infty$.

Following [Iy3], assume that the Abelian category $\mathscr{A}=\mod
\Lambda$ and $\mathscr{B}={^{\bot}\Lambda}$. Then the categories
${^{\bot_{n}}\Lambda}\bigcap{^{\bot}\Lambda}=\mathscr{B}$ and
${\Lambda^{\bot_{n}}}\bigcap{^{\bot}\Lambda}=\mathscr{B}$ for any
$n\geq 1$. Denote by $\underline{\mathscr{B}}$ (resp.
$\overline{\mathscr{B}}$) the stable category $\mathscr{B}$ modulos
relative projectives (resp. injectives) in $\mathscr{B}$. We remark
that $\mathscr{B}$ forms a Frobenius category, so the relative
projectives in $\mathscr{B}$ coincide with relative injectives in
$\mathscr{B}$, and that we have
$\underline{\mathscr{B}}=\overline{\mathscr{B}}$.

We use $\Proj (\mod \overline{\mathscr{B}})$ to denote the
subcategory of $\mod \overline{\mathscr{B}}$ consisting of
projective objects. The functors
$\tau:\underline{\mathscr{B}}\rightarrow \overline{\mathscr{B}}$ and
$\tau^{-}:\overline{\mathscr{B}}\rightarrow \underline{\mathscr{B}}$
are quasi-inverse equivalences, where $\tau=F^{-}\circ G$ and
$\tau^{-}=G^{-}\circ F$ with $F:\overline{\mathscr{B}}\rightarrow
\Proj (\mod \overline{\mathscr{B}})$ via $X\rightarrow
\underline{\mathscr{B}}(\ ,X)$ and $G:
\underline{\mathscr{B}}\rightarrow \Proj (\mod
\overline{\mathscr{B}})$ via $X\rightarrow \mathbb{D}{\rm
Ext}_{\mathscr{A}}^{1}(X,\ )$. Let $\Omega$ and $\Omega^{-}$ be the
relative syzygy and cosyzygy functors in ${^{\bot}\Lambda}$.

\vspace{0.2cm}

{\bf Lemma 4.4} (1) ([Iy3, Corollary 2.3.2]) {\it The functors
$\tau$ and $\tau^{-}$ give mutually-inverse equivalences
$\tau:\underline{\mathscr{B}}\rightarrow \underline{\mathscr{B}}$
and $\tau^{-}:\underline{\mathscr{B}}\rightarrow
\underline{\mathscr{B}}$.

(2) {\rm ([Iy3, Theorem 2.3.1])} There exist functorial isomorphisms
$\underline{{\rm Hom}}_{\mathscr{B}}(Y,\tau X)\cong \mathbb{D}{\rm
Ext}^{1}_{\mathscr{A}}(X,Y)$ \linebreak $\cong \underline{{\rm
Hom}}_{\mathscr{B}}(\tau^{-}Y,X)$ for any $X,Y\in
\underline{\mathscr{B}}$.}

\vspace{0.2cm}

The following result is a generalization of [EH, Lemma 3.2].

\vspace{0.2cm}

{\bf Lemma 4.5} {\it For any $M,N\in {^{\bot}\Lambda}$ and $i\geq1$,
${\rm Ext}^{i}_{\mathscr{A}}(M,N)\cong \mathbb{D}{\rm
Ext}^{1}_{\mathscr{A}}(N,\Omega^{i}\tau M)$.}

\vspace{0.2cm}

{\it Proof.} Note that the functors
$\Omega:\underline{\mathscr{B}}\rightarrow\underline{\mathscr{B}}$
and
$\Omega^{-}:\underline{\mathscr{B}}\rightarrow\underline{\mathscr{B}}$
are mutually-inverse equivalences by [AR2, Proposition 3.1]. Then by
Lemma 4.4(2), we have that ${\rm Ext}_{\mathscr{A}}^{i}(M,N)\cong
\underline{{\rm Hom}}_{\mathscr{B}}(\Omega^{i}M,N)\cong
\underline{{\rm Hom}}_{\mathscr{B}}(\Omega M,\Omega^{-i+1}N)\cong
{\rm Ext}_{\mathscr{A}}^{1}(M,\Omega^{-i+1}N)\cong
\mathbb{D}\underline{{\rm
Hom}}_{\mathscr{B}}(\tau^{-1}\Omega^{-i+1}N,M)$ \linebreak $\cong
\mathbb{D}\underline{{\rm Hom}}_{\mathscr{B}}(\Omega^{-i+1}N,\tau
M)\cong \mathbb{D}\underline{{\rm
Hom}}_{\mathscr{B}}(\Omega^{1}N,\Omega^{i}\tau M)\cong
\mathbb{D}{\rm Ext}_{\mathscr{A}}^{1}(N,\Omega^{i}\tau M)$.
$\hfill{\square}$

\vspace{0.2cm}

Let $M$ be a module in $\mod \Lambda$. For a positive integer $n$,
$M\in \mod \Lambda$ is called $\Omega^{n}\tau$-{\it periodic} if
there exists a positive integer $t$, such that $(\Omega^{n}\tau)^tM
\cong M$. Recall from [EH] that $M$ is called a {\it maximal}
$n$-{\it orthogonal module} if $\add_{\Lambda}M$ is a maximal
$n$-orthogonal subcategory of $\mod \Lambda$. The following result
generalizes [EH, Theorem 3.1]. It should be pointed out that this
result can be induced by [Iy3, Theorem 2.5.1(1)].

 \vspace{0.2cm}

{\bf Theorem 4.6} {\it Let $X \in {^{\bot}\Lambda}$ be a maximal
$n$-orthogonal module. If $Y$ is a direct summand of $X$ then so is
$\Omega^{n}\tau Y$ for some $n\geq 1$. Hence the non-projective
direct summand of $X$ is $\Omega^{n}\tau$-periodic.}

\vspace{0.2cm}

{\it Proof.} By Lemma 4.4(1), $\tau X\in {^{\bot}\Lambda}$ and so
$\Omega ^{n}\tau X\in {^{\bot}\Lambda}$. Since $X$ is maximal
$n$-orthogonal, for any $1\leq i\leq n$, we have that ${\rm
Ext}_{\mathscr{A}}^{i}(X,\Omega^{n}\tau X)\cong \underline{{\rm
Hom}}_{\mathscr{B}}(\Omega^{i}X,\Omega^{n}\tau X)\cong
\underline{{\rm Hom}}_{\mathscr{B}}(\Omega^{1}X,\Omega^{-i+n+1}\tau
X)\cong {\rm Ext}_{\mathscr{A}}^{1}(X,\Omega^{-i+n+1}\tau X)\cong
\mathbb{D}{\rm Ext}_{\mathscr{A}}^{-i+n+1}(X, X)=0$ by Lemma 4.5.
Also since $X$ is maximal $n$-orthogonal, $\Omega^{n}\tau X\in
\add_{\Lambda}X$. Notice that both $\Omega$ and $\tau$ are
equivalences in $^{\bot}\Lambda$, it follows that if $Y$ is a direct
summand of $X$, then $\Omega^{n}\tau Y$ is also a direct summand of
$X$. Since $X$ has only finitely many indecomposable direct
summands, some power of $\Omega^{n}\tau$ is identity on the
non-projective direct summands of $X$. $\hfill{\square}$

\vspace{0.5cm}

\centerline{\large\bf 5. The Complexity of Modules}

\vspace{0.2cm}

In this section, both $\Lambda$ and $\Gamma$ are finite-dimensional
$K$-algebras over a field $K$. We will study the relation between
the complexity of modules and the existence of maximal
$n$-orthogonal subcategories of $\mod \Lambda \otimes _{K} \Gamma$.

From the following result, we can construct algebras with infinite
global dimension admitting no maximal $n$-orthogonal subcategories
for any $n \geq 1$.

\vspace{0.2cm}

{\bf Proposition 5.1} {\it If $\id _{\Lambda}\Lambda=n \geq 1$ and
${\Lambda \otimes_{K}\Gamma}$ admits a maximal $j$-orthogonal
subcategory of $\mod {\Lambda \otimes_{K}\Gamma}$ for some $j \geq
n$, then ${\rm Hom}_{\Gamma}(\mathbb{D}\Gamma ^{op},\Gamma)=0$.}

\vspace{0.2cm}

{\it Proof.} Suppose that $\mathscr{C}$ is a maximal $j$-orthogonal
subcategory of $\mod {\Lambda\otimes_{K}\Gamma}$ for some $j\geq n$.
Then both ${\Lambda\otimes_{K}\Gamma}$ and ${\mathbb{D}\Lambda
^{op}\otimes_{K}\mathbb{D}\Gamma ^{op}}$ are in $\mathscr{C}$. Thus
by [CE, Chapter XI, Theorem 3.1], $0={\rm
Ext}_{{\Lambda\otimes_{K}\Gamma}}^{i}({\mathbb{D}\Lambda ^{op}
\otimes_{K}\mathbb{D}\Gamma ^{op}},{\Lambda\otimes_{K}\Gamma})=
\bigoplus_{r+s=i}{\rm Ext}_{\Lambda}^{r}(\mathbb{D}\Lambda
^{op},\Lambda)\otimes_{K}{\rm Ext}_{\Gamma}^{s}(\mathbb{D}\Gamma
^{op},\Gamma)$ for any $1\leq i\leq n$. So ${\rm
Ext}_{\Lambda}^{i}(\mathbb{D}\Lambda ^{op},\Lambda)\otimes_{K}{\rm
Hom}_{\Gamma}(\mathbb{D}\Gamma ^{op},\Gamma)=0$ for any $1\leq i\leq
n$.

If $0\neq{\rm Hom}_{\Gamma}(\mathbb{D}\Gamma ^{op},\Gamma)$, then
${\rm Hom}_{\Gamma}(\mathbb{D}\Gamma ^{op},\Gamma)\cong K^{m}$ as
$K$-vector spaces for some $m \geq 1$. Thus $[{\rm
Ext}_{\Lambda}^{i}(\mathbb{D}\Lambda ^{op},\Lambda)]^m \cong {\rm
Ext}_{\Lambda}^{i}(\mathbb{D}\Lambda ^{op},\Lambda)\otimes_{K}K^m
\cong {\rm Ext}_{\Lambda}^{i}(\mathbb{D}\Lambda
^{op},\Lambda)\otimes_{K}{\rm Hom}_{\Gamma}(\mathbb{D}\Gamma
^{op},\Gamma)=0$ and ${\rm Ext}_{\Lambda}^{i}(\mathbb{D}\Lambda
^{op},\Lambda)=0$ for any $1\leq i\leq n$. Since $\mathbb{D}\Lambda
^{op}$ is an injective cogenerator for $\mod \Lambda$ and $\id
_{\Lambda}\Lambda=n$, it is easy to see that $\Lambda$ is
selfinjective by applying the functor ${\rm
Hom}_{\Lambda}(\mathbb{D}\Lambda ^{op},-)$ to the minimal injective
resolution of $_{\Lambda}\Lambda$, which is a contradiction.
\hfill{$\square$}

\vspace{0.2cm}

Assume that $\id _{\Lambda}\Lambda=1$. If there exists a torsionless
injective module in $\mod \Gamma$ (that is, there exists an
injective module $Q\in \mod \Gamma$ such that the canonical
evaluation homomorphism $Q \rightarrow {\rm Hom}_{\Gamma}({\rm
Hom}_{\Gamma}(Q, \Gamma), \Gamma)$ is monomorphic), especially, if
there exists a projective-injective module in $\mod \Gamma$ (for
example, $\Gamma$ is 1-Gorenstein), then
${\Lambda\otimes_{K}\Gamma}$ admits no maximal $n$-orthogonal
subcategories of $\mod {\Lambda\otimes_{K}\Gamma}$ for any $n \geq
1$ by Proposition 5.1.

\vspace{0.2cm}

{\bf Definition 5.2} ([EH]) Let
$$\cdots\rightarrow P_{n}\rightarrow\cdots\rightarrow
P_{1}\rightarrow P_{0}\rightarrow M\rightarrow 0$$ be a minimal
projective resolution of a module $M \in \mod \Lambda$. The {\it
complexity} of $M$ is defined as $\cx(M)= \inf \{ b\geq 0\ |$ there
exists a $c>0$ such that $\dim _{K}P_{n}\leq cn^{b-1}$ for all $n\}$
if it exists, otherwise $\cx(M)=\infty$.

\vspace{0.2cm}

It is easy to see that $\cx(M)=0$ implies $M$ is of finite
projective dimension, and $\cx(M)\leq1$ if and only if the
dimensions of $P_{n}$ are bounded. Erdmann and Holm proved in [EH,
Theorem 1.1] that if $\Lambda$ is selfinjective and there exists a
maximal $n$-orthogonal module in $\mod \Lambda$ for some $n \geq 1$,
then all modules in $\mod \Lambda$ have complexity at most 1. At the
end of [EH], Erdmann and Holm posed a question: Whether there can
exist maximal $n$-orthogonal modules for non-selfinjective algebras
$\Lambda$ for which not all modules in $\mod \Lambda$ are of
complexity at most 1? In the following, we will give some properties
of the complexity of the tensor product of modules. Then we
construct a class of algebras $\Lambda$ with $\id
_{\Lambda}\Lambda=n(\geq 1)$, such that not all modules in $\mod
\Lambda$ are of complexity at most 1, but $\Lambda$ admits no
maximal $j$-orthogonal subcategories of $\mod \Lambda$ for any
$j\geq n$.

\vspace{0.2cm}

{\bf Lemma 5.3} {\it Let $$\cdots\rightarrow
P_{i}\stackrel{f_{i}}{\rightarrow}\cdots\rightarrow
P_{1}\stackrel{f_{1}}{\rightarrow}
P_{0}\stackrel{f_{0}}{\rightarrow} M\rightarrow 0$$ and
$$\cdots\rightarrow Q_{j}\stackrel{g_{j}}
{\rightarrow}\cdots\rightarrow Q_{1}\stackrel{g_{1}}{\rightarrow}
Q_{0}\stackrel{g_{0}}\rightarrow N\rightarrow 0$$ be minimal
projective resolutions of $M \in \mod \Lambda$ and $N \in \mod
\Gamma$, respectively. Then the following is a minimal projective
resolution of $M\otimes_{K}N$ as a
$\Lambda\otimes_{K}\Gamma$-module:
$$\cdots\rightarrow R_{n}\rightarrow \cdots \rightarrow R_{1}
\stackrel{(1_{P_{0}}\otimes_{K} g_{1})\bigoplus (f_{1} \otimes
_{K}1_{Q_{0}})}{\longrightarrow} R_{0}\stackrel{f_{0}\otimes_{K}
g_{0}}{\longrightarrow} M\otimes_{K}N\rightarrow 0 \eqno{(2)}$$
where $R_{n}=\bigoplus _{i+j=n}P_{i}\otimes_{K}Q_{j}$ for any $n
\geq 0$.}

\vspace{0.2cm}

{\it Proof.} It is well known that the sequence (2) is a projective
resolution of $M\otimes_{K}N$ {\it as a}
$\Lambda\otimes_{K}\Gamma$-module. So it suffices to prove the
minimality.

It is straightforward to verify that $\Ker f_{0}\otimes_{K} g_{0}
\cong (\Ker f_{0}\otimes_{K} Q_{0})+(P_{0}\otimes_{K} \Ker g_{0})$.
On the other hand, note that $\Ker f_{0}\subseteq J(\Lambda)P_{0}$,
$\Ker g_{0}\subseteq J(\Gamma)Q_{0}$ and the nilpotent ideal
$(J(\Lambda)\otimes _{K} \Gamma)+(\Lambda\otimes_{K}
J(\Gamma))\subseteq J(\Lambda\otimes_{K} \Gamma)$. So we have that
$(\Ker f_{0}\otimes_K Q_{0})+(P_{0}\otimes_K \Ker g_{0})\subseteq
J(P_{0}\otimes_K Q_{0})$ and hence $f_{0}\otimes_{K} g_{0}$ is
minimal. By using an argument similar to above, we get the desired
assertion. $\hfill{\square}$

\vspace{0.2cm}

{\bf Proposition 5.4} {\it $\max\{\cx(M),\cx(N)\}\leq
\cx(M\otimes_{K}N)\leq \cx(M)+\cx(N)$ for any $M \in \mod \Lambda$
and $N \in \mod \Gamma$.}

\vspace{0.2cm}

{\it Proof.} By Lemma 5.3, we have that $\dim_{K}
R_{n}=\sum_{i+j=n}\dim_{K} P_{i}\dim_{K} Q_{j}\geq$ \linebreak
$\dim_{K} P_{n}\dim_{K} Q_{0}\geq \dim_{K} P_{n}$. Similarly,
$\dim_{K} R_{n}\geq \dim_{K} Q_{n}$. So $\max\{\cx(M), \cx(N)\}\leq
\cx(M\otimes_{K}N)$.

To prove the second inequality, without loss of generality, assume
that both $\cx(M)$ and $\cx(N)$ are finite. Then there exist $c, c'
>0$ such that $\dim_{K} P_{i}\leq
ci^{\cx(M)-1}$ and $\dim_{K} Q_{j}\leq c'j^{\cx(N)-1}$ for any $i,
j\geq 0$. Thus by Lemma 5.3, we have that $\dim_{K}
R_{n}=\sum_{i+j=n}{\dim_{K} P_{i}\dim_{K} Q_{j}}\leq
\sum_{i+j=n}(ci^{\cx(M)-1})(c'j^{\cx(N)-1})\leq
(n+1)(cc')n^{\cx(M)+\cx(N)-2}\leq (2cc')n^{\cx(M)+\cx(N)-1}$, which
implies $\cx(M\otimes_{K}N)\leq \cx(M)+\cx(N)$. $\hfill{\square}$

\vspace{0.2cm}

Now we are in a position to give the following example.

\vspace{0.2cm}

{\bf Example 5.5} Let $\id _{\Lambda}\Lambda=n\geq 1$ and $\Gamma$
be selfinjective of infinite representation type with
$J(\Gamma)^3=0$. Then we have

(1) $\id _{\Lambda\otimes_{K}\Gamma}\Lambda\otimes_{K}\Gamma =n$.

(2)  $\sup \{ \cx(X) | X\in \mod {\Lambda\otimes_{K}\Gamma} \} \geq
2$.

(3) ${\Lambda\otimes_{K}\Gamma}$ admits no maximal $j$-orthogonal
subcategories of $\mod {\Lambda\otimes_{K}\Gamma}$ for any $j\geq
n$.

\vspace{0.2cm}

{\it Proof.} Because $\Gamma$ is selfinjective, ${\rm
Hom}_{\Gamma}(\mathbb{D}\Gamma ^{op}, \Gamma)\neq 0$ and then the
assertion (3) follows from Proposition 5.1. It is well known that
$\max\{\id _{\Lambda}\Lambda,\id _{\Gamma}\Gamma \}\leq \id
_{\Lambda\otimes_{K}\Gamma}\Lambda\otimes_{K}\Gamma\leq \id
_{\Lambda}\Lambda +\id _{\Gamma}\Gamma$ (cf. [AR2, Proposition
2.2]). By assumption, $\id _{\Lambda}\Lambda =n$ and $\Gamma$ is
selfinjective, so $\id
_{\Lambda\otimes_{K}\Gamma}\Lambda\otimes_{K}\Gamma =n$ and the
assertion (1) follows. By [GW, Theorem 6.1], $\sup \{ \cx(N) | N\in
\mod \Gamma \} \geq 2$. Then it follows from Proposition 5.4 that
the assertion (2) holds true. $\hfill{\square}$

\vspace{0.5cm}

{\bf Acknowledgements} The research was partially supported by the
Specialized Research Fund for the Doctoral Program of Higher
Education (Grant No. 20060284002), NSFC (Grant No.10771095) and NSF
of Jiangsu Province of China (Grant No. BK2007517). The authors
thank Prof. Osamu Iyama and the referee for the useful suggestions.


\vspace{0.5cm}

\begin{thebibliography}{101}

\bibitem[AR1]{A1} M. Auslander and I. Reiten, {\it Applications of
contravariantly finite subcategories}. Adv. Math. {\bf 86} (1991),
111--152.

\bibitem[AR2]{A2} M. Auslander and I. Reiten, {\it Cohen-Macaulay and
Gorenstein algebras}. In: Representation theory of finite groups and
finite-dimensional algebras, Bielefeld, 1991, edited by G.O. Michler
and C.M. Ringel, Prog. Math., {\bf 95}, Birkh\"auser, Basel, 1991,
pp. 221--245.

\bibitem[BR]{A3} A. Beligiannis and I. Reiten, {\it Homological and
homotopical aspects of torsion theories}, Memoirs Amer. Math. Soc.,
{\bf 188}, Amer. Math. Soc., Providence, RI. 2007.

\bibitem[Bj]{A4} J. E. Bj\"ork, {\it The Auslander condition on
Noetherian rings}. In: S\'{e}minaire d'Alg\`{e}bre Paul Dubreil et
Marie-Paul Malliavin, 39\`{e}me Ann\'{e}e, Paris, 1987/1988, Lect.
Notes in Math., {\bf 1404}, Springer-Verlag, Berlin, 1989,
pp.137--173.

\bibitem[CE]{A5} H. Cartan and S. Eilenberg, Homological Algebra.
Reprint of the 1956 original, Princeton Landmarks in Math.,
Princeton Univ. Press, Princeton, 1999.

\bibitem[EH]{A6} K. Erdmann and T. Holm, {\it Maximal n-orthogonal
modules for selfinjective algebras}. Proc. Amer. Math. Soc. {\bf
136} (2008), 3069--3078.

\bibitem[FGR]{A7} R. M. Fossum, P.A. Griffith and I. Reiten, Trivial
Extensions of Abelian Categories. Lect. Notes in Math. {\bf 456},
Springer-Verlag, Berlin, 1975.

\bibitem[GLS]{A8} C. Geiss, B. Leclerc and J. Schr\"oer, {\it Rigid
modules over preprojective algebras}. Invent. Math. {\bf 165}
(2006), 589-632.

\bibitem[GW]{A9} J. Y. Guo and Q. X. Wu, {\it Selfinjective Koszul
algebras of finite complexity}. Preprint.

\bibitem[HRS]{A10} D. Happel, I. Reiten and S.O. Smal$\phi$,
{\it Tilting in Abelian categories and Quasitilted Algebras}.
Memoirs Amer. Math. Soc., {\bf 575}, Amer. Math. Soc., Providence,
RI, 1996.

\bibitem[Ho]{A11} M. Hoshino, {\it Syzygies and Gorenstein rings}.
Arch. Math. (Basel) {\bf 55} (1990), 355--360.

\bibitem[IS]{A12} Y. Iwanaga, H. Sato, {\it On Auslander's}
$n$-{\it Gorenstein rings}. J. Pure Appl. Algebra {\bf 106} (1996),
61--76.

\bibitem[Iy1]{A13} O. Iyama, {\it Symmetry and duality on n-Gorenstein
rings}. J. Algebra {\bf 269} (2003), 528--535.

\bibitem[Iy2]{A14} O. Iyama, {\it Higher-dimensional Auslander-Reiten
theory on maximal orthogonal subcategories}. Adv. Math. {\bf 210}
(2007), 22--50.

\bibitem[Iy3]{A15} O. Iyama, {\it Auslander correspondence}. Adv. Math. {\bf 210}
(2007), 51--82.

\bibitem[Iy4]{A16} O. Iyama, {\it Cluster tilting for higher Auslander algebras}.
Preprint is available at: arXiv:0809.4897.

\bibitem[M]{A17} Y. Miyashita, {\it Tilting modules of finite
projective dimension}. Math. Z. {\bf 193} (1986), 113--146.

\end{thebibliography}
\end{document}